\documentstyle[12pt]{article}
\newfam\msbfam
\font\tenmsb=msbm10    \textfont\msbfam=\tenmsb \font\sevenmsb=msbm7
\scriptfont\msbfam=\sevenmsb \font\fivemsb=msbm5
\scriptscriptfont\msbfam=\fivemsb
\def\Bbb{\fam\msbfam \tenmsb}

\def\rr{{\Bbb R}}
\def\rz{{{\rr}^n}}

\def\fz{\infty}
\def\az{\alpha}
\def\supp{{\rm{\ supp\ }}}

\def\ez{\epsilon}
\def\bz{\beta}
\def\pz{\partial}

\def\sz{\sigma}
\def\vz{\varphi}
\def\dz{\delta}
\def\gz{\gamma}

\def\lz{\lambda}

\def\supp{{\rm supp}}

\def\wt{\widetilde}
\def\wz{\omega}
\def\l{\left}
\def\r{\right}

\def\dsum{\displaystyle\sum}

\def\dint{\displaystyle\int}
\def\dfrac{\displaystyle\frac}
\def\dsup{\displaystyle\sup}

\def\dinf{\displaystyle\inf}

\newtheorem{thm}{\hskip\parindent Theorem}

\newtheorem{lem}{\hskip\parindent Lemma}

\newtheorem{prop}{\hskip\parindent Proposition}

\newtheorem{cor}{\hskip\parindent Corollary}

\begin{document}

\baselineskip=15pt
\renewcommand{\arraystretch}{2}
\arraycolsep=1.2pt

\title{ Weighted norm inequalities for pseudo-differential
operators with smooth symbols and their commutators} {\footnotetext{ \hspace{-0.65
cm} 2000 Mathematics Subject  Classification:  42B20, 35S05, 47G30.\\
The  research was supported  by the NNSF (10401002) of China.\\}

\author{ Lin Tang}
\date{}
\maketitle

{\bf Abstract}\quad We obtain weighted $L^p$ inequalities for
pseudo-differential operators with smooth symbols and their
commutators by using  a class of new  weight functions which
include Muckenhoupt weight functions. Our results  improve
essentially some well-known results.

\bigskip

\begin{center}{\bf 1. Introduction }\end{center}
Let $m$ be real number. Following \cite{t}, a symbol in
$S_{1,\dz}^m$ is a smooth function $\sz(x,\xi)$ defined on
$\rz\times\rz$ such that for all multi-indices $\az$ and $\bz$ the
following estimate holds:
$$|D_x^\az D_\xi^\bz\sz(x,\xi)|\le
C_{\az,\bz}(1+|\xi|)^{m-|\bz|+\dz|\az|},$$ where $ C_{\az,\bz}>0$
is independent of $x$ and $\xi$. A symbol in $S_{1,\dz}^{-\fz}$ is
one which satisfies the above estimates for each real number $m$.

The operator $T$ given by
$$Tf(x)=\dint_\rz\sz(x,\xi)e^{2\pi ix\cdot \xi}\hat f(\xi)\,d\xi$$
is called a pseudo-differential operator with symbol
$\sz(x,\xi)\in S_{1,\dz}^m$, where $f$ is a Schwartz function and
$\hat f$ denotes the Fourier transform of $f$. As usual,
$L^m_{1,\dz}$ will denote the class of pseudo-differential
operators with symbols in $S_{1,\dz}^m$.

Miller \cite{m} showed the boundedness of $L^0_{1,0}$
pseudo-differential operators on weighted $L^p(1<p<\fz)$ spaces
whenever the weight function belongs to Muckenhoupt's class $A_p$.
In particular, A. Laptev \cite{l} proved that any  $L^0_{1,0}$
pseudo-differential operator is a standard Calder\'on-Zygmund
operator.

Our purpose is to improve  the  above  results of Miller. To state
the main results, let us first introduce some notations.

In this paper,  $Q(x,t)$ denotes the cube centered at $x$ and of the
sidelength $t$. Similarly, given $Q=Q(x,t)$ and $\lz>0$, we will
write $\lz Q$ for the $\lz$-dilate cube, which is the cube with the
same center $x$ and with sidelength $\lz t$. Given a Lebesgue
measurable set $E$ and a weight $\wz$, $|E|$ will denote the
Lebesgue measure of $E$ and $\wz(E)=\int_E\wz dx$.
$\|f\|_{L^p(\wz)}$ will denote $(\int_\rz |f(y)|^p\wz(y)dy)^{1/p}$
for $0< p<\fz$. $C$ denotes the constants that are independent of
the main parameters involved but whose value may differ from line to
line. For a measurable set $E$, denote by $\chi_E$ the
characteristic function of $E$. By $A\sim B$, we mean that there
exists a constant $C>1$ such that $1/C\le A/B\le C$.

Throughout this paper, we let $\vz(t)=(1+t)^{\az_0}$ for $\az_0>0$
and $t\ge 0$.

 A weight will always mean a positive function which is
locally integrable. We say that a weight $\wz$ belongs to the class
$A_p(\vz)$ for $1<p<\fz$, if there is a constant $C$ such that for
all cubes  $Q=Q(x,r)$ with center $x$ and sidelength $r$
$$\l(\dfrac 1{\vz(|Q|)|Q|}\dint_Q\wz(y)\,dy\r)
\l(\dfrac 1{\vz(|Q|)|Q|}\dint_Q\wz^{-\frac 1{p
-1}}(y)\,dy\r)^{p-1}\le C.$$
 We also
say that a  nonnegative function $\wz$ satisfies the $A_1(\vz)$
condition if there exists a constant $C$ for all cubes $Q$
$$M_\vz(\wz)(x)\le C \wz(x), \ a.e.\ x\in\rz.$$
 where
$$M_\vz f(x)=\dsup_{x\in Q}\dfrac 1{\vz(|Q|)|Q|}\dint_Q|f(y)|\,dy.$$
We also defined the Hardy-Littlewood maximal operator $M$ by
$$M f(x)=\dsup_{x\in Q}\dfrac 1{|Q|}\dint_Q|f(y)|\,dy.$$
Obviously, $f(x)\le M_\vz f(x)\le M f(x)$ a.e. $x\in\rz$ and the
function $M_\vz f(x)$ is lower semi-continuous.

   Since $\vz(|Q|)\ge 1$, so
$A_p(\rz)\subset A_p(\vz)$ for $1\le p<\fz$, where $A_p(\rz)$ denote
the classical Muckenhoupt weights; see \cite{gr}. It is well known
that if $\wz\in A_\fz(\rz)=\bigcup_{p\ge 1}A_p(\rz)$, then $\wz(x)
dx$ be  a doubling measure, that is, there exist a constant $C>0$
for any cube $Q$ such that
$$\wz(2Q)\le C\wz(Q).$$
 From the definition of $A_p(\vz)$ and (iv) of
Lemma 2.1 in Section 2, it is easy to see that if $\wz\in
A_p(\vz)$, then $\wz(x)dx$ may be not a doubling measure. In fact,
let $0\le\gz\le n\az_0$, it is easy to check that
$\wz(x)=(1+|x|)^{-(n+\gz)}\not\in A_\fz(\rz)$ and $\wz(x)dx$ is
not a doubling measure, but $\wz(x)=(1+|x|)^{-(n+\gz)}\in
A_1(\vz)$.

Now let us state our  main results as follows.
\begin{thm}\label{t1.1.}\hspace{-0.1cm}{\rm\bf 1.1.}\quad
Suppose $T\in L_{1,0}^0$. Let $1<p<\fz$ and $\wz\in A_{p}(\vz)$,
then there exists a constant $C>0$ such that
$$\|Tf\|_{L^p(\wz)}\le C\|f\|_{L^p(\wz)}.$$
\end{thm}
As a consequence of Theorem 1.1, we have the following result.
\begin{cor}\label{c1.1.}\hspace{-0.1cm}{\rm\bf 1.1.}\quad
Suppose $T\in L_{1,0}^0$. Let $1<p<\fz$ and
$\wz(x)=(1+|x|)^{\gz_1}  |x|^{\gz_2}$ with $\gz_1\in\rr$ and
$-n<\gz_2<n(p-1)$. Then there exists a constant $C>0$ such that
$$\|Tf\|_{L^p(\wz)}\le C\|f\|_{L^p(\wz)}.$$
\end{cor}

Let $b\in BMO$(bounded mean oscillation function) defined in
\cite{jn}. We also consider the commutator of
Coifman-Rochberg-Weiss $[b,T]$ defined by
$$[b,T]f(x)=b(x)Tf(x)-T(bf)(x).$$
Similar to Theorem 1.1, we have
\begin{thm}\label{t1.2.}\hspace{-0.1cm}{\rm\bf 1.2.}\quad
Suppose $T\in L_{1,0}^0$. Let $b\in BMO$, $1< p<\fz$ and $\wz\in
A_{p}(\vz)$, then  there exists a constant $C>0$ such that
$$\|[b,T]f\|_{L^p(\wz)}\le C\|b\|_{BMO}\|f\|_{L^p(\wz)}.$$
\end{thm}
As a consequence of Theorem 1.2, we have the following result.
\begin{cor}\label{c1.2.}\hspace{-0.1cm}{\rm\bf 1.2.}\quad
Suppose $T\in L_{1,0}^0$. Let $b\in BMO$, $1<p<\fz$ and
$\wz(x)=(1+|x|)^{\gz_1}|x|^{\gz_2}$ with $\gz_1\in\rr$ and
$-n<\gz_2<n(p-1)$. Then there exists a constant $C>0$ such that
$$\|[b,T]f\|_{L^p(\wz)}\le C\|b\|_{BMO}\|f\|_{L^p(\wz)}.$$
\end{cor}
We remark that the above results also hold if $T\in L_{1,\dz}^0$
with $0<\dz<1$.

The organization of the paper is as follows: We study some
elementary properties of the new weight functions in Section 2. We
proved Theorem 1.1 and Corollary 1.1 in Section 3. Theorem 1.2 and
Corollary 1.2 are proved, and  the weighted weak (1,1) type
inequality is obtained in Section 4.

Remark: We will consider the similar results for
bilinear pseudo-differential
operators in the forthcoming paper.

\begin{center} {\bf 2. Some properties of $A_p(\vz)$}\end{center}

Similar to the classical Muckenhoupt weights, we give some
properties for weights $\wz\in A_\fz(\vz)=\bigcup_{p\ge
1}A_p(\vz)$.
\begin{lem}\label{l2.1.}\hspace{-0.1cm}{\rm\bf 2.1.}\quad For any cube
$Q\subset \rz$, then
\begin{enumerate}
\item[(i)]If $ 1\le p_1<p_2<\fz$, then $A_{p_1}(\vz)\subset
A_{p_2}(\vz)$. \item[(ii)] $\wz\in A_p(\vz)$ if and only if
$\wz^{-\frac 1{p-1}}\in A_{p'}(\vz)$, where $1/p+1/p'=1.$
\item[(iii)] If $\wz_1,\ \wz_2\in A_p(\vz),\ p\ge1$, then
$\wz_1^\az\wz_2^{1-\az}\in A_p(\vz)$ for any $0<\az<1$.
\item[(iv)] If $\wz\in A_p$ for $1\le p<\fz$, then
$$\dfrac 1{\vz(|Q|)|Q|}\dint_Q|f(y)|dy\le C\l(\dfrac
1{\wz(5Q)}\dint_Q|f|^p\wz(y)dy\r)^{1/p}.$$  In particular, let
$f=\chi_{E}$ for any measurable set $E\subset Q$,
$$\dfrac {|E|}{\vz(|Q|)|Q|}\le C\l(\dfrac
{\wz(E)}{\wz(5Q)}\r)^{1/p}.$$
\end{enumerate}
\end{lem}
Proof:\quad (i), (ii) and (iii) are obvious. We only prove (iv). In
fact,
$$\begin{array}{cl}
\dfrac 1{\vz(|Q|)|Q|}&\dint_Q|f(y)|dy\\&=\dfrac
1{\vz(|Q|)|Q|}\dint_Q|f(y)|\wz^{\frac 1p}(y)\wz^{-\frac 1p}(y)dy\\
&\le \l(\dfrac 1{\vz(|Q|)|Q|}\dint_Q|f(y)|^p\wz(y)dy\r)^{\frac
1p}\l(\dfrac
1{\vz(|Q|)|Q|}\dint_Q\wz^{-\frac 1{p-1}}(y)dy\r)^{\frac {p-1}p}\\
&\le C\l(\dfrac 1{\vz(|Q|)|Q|}\dint_Q|f(y)|^p\wz(y)dy\r)^{\frac
1p}\\
&\qquad\qquad\times\l(\dfrac
1{\vz(|5Q|)|5Q|}\dint_{5Q}\wz^{-\frac 1{p-1}}(y)dy\r)^{\frac {p-1}p}\\
&\le C\l(\dfrac 1{\vz(|Q|)|Q|}\dint_Q|f(y)|^p\wz(y)dy\r)^{\frac
1p}\l(\dfrac
1{\vz(|5Q|)|5Q|}\dint_{5Q}\wz(y)dy\r)^{-\frac 1p}\\
&\le C\l(\dfrac 1{\wz(5Q)}\dint_Q|f(y)|^p\wz(y)dy\r)^{1/p}.
\end{array}$$
Thus, (iv) is proved.

Obviously, it is easy to see that $\cal S$ (the set of all Schwartz
functions) is dense in $L^p(\wz)$ for $\wz\in A_\fz(\vz)$ and $1\le
p<\fz$. Hence, we always assume $f\in  \cal S$ if  $f\in L^p(\wz)$
for $1\le p<\fz$.

 Next, we give a result about the operator $M_{\wz}$ defined by
$$M_{\wz}(f)(x)=\dsup_{x\in Q}\dfrac
1{\wz(5Q)}\dint_Q|f(x)|\wz(x)dx.$$
\begin{lem}\label{l2.2.}\hspace{-0.1cm}{\rm\bf 2.2.}\quad
Let $\wz\in A_\fz(\vz )$, then
$$\wz(\{x\in\rz:\ M_{\wz}f(x)>\lz\})\le \dfrac
C\lz\|f\|_{L^1(\wz)},\ \forall\lz>0,\ \forall f\in
L^1(\wz),\eqno(2.1)$$ and for $1<p<\fz$
$$\|M_{\wz}f\|_{L^p(\wz)}\le C_p\|f\|_{L^p(\wz)}.\eqno(2.2)$$
\end{lem}
 Proof:\quad  We  set $x\in
E_\lz=\{y\in\rz:\ M_{\wz}f(y)>\lz\}$ with any $\lz>0$, then, there
exists a cube $Q_x\ni x$ such that
$$\dfrac
1{\wz(5Q_x)}\dint_{Q_x}|f(x)|\wz(x)dx>\lz.$$ Thus, $\{Q_x\}_{x\in
E_\lz}$ covers $E_\lz$. By Vitali lemma, there exists a class
disjoint cubes $\{Q_{xj}\} $ such that $\bigcup Q_{xj}\subset
E_\lz\subset \bigcup 5Q_{xj}$ and
$$\wz(E_\lz)\le \dsum_j\wz(5Q_{xj})\le \dfrac
C\lz\dsum_j\dint_{Q_{xj}}|f(y)|\wz(y)dy\le \dfrac
C\lz\|f\|_{L^1(\wz)}.$$ Thus, (2.1) is proved. By interpolation
from (2.1) and the boundedness of $M_{\wz}$ in $L^\fz(\wz)$, we
obtain (2.2). The proof is finished.

From (iii) of Lemma 2.1, we know that
$$M_\vz f(x)\le C(M_{\wz}(|f|^p)(x))^{1/p},\quad x\in\rz.$$
From this and using Lemma 2.2, we can get the following result.
\begin{lem}\label{l2.3.}\hspace{-0.1cm}{\rm\bf 2.3.}\quad
Let $1\le p_1 <\fz$ and suppose that $\wz\in A_{p_1}(\vz)$. If
$p_1<p<\fz$, then the equality
$$\dint_\rz|M_\vz f(x)|^{p}\wz(x)dx\le
C_{p}\dint_\rz|f(x)|^{p}\wz(x)dx.$$ Furthermore, let $1\le p<\fz$,
$\wz\in A_p(\vz)$   if and only if
$$\wz(\{x\in\rz:\ M_{\vz}f(x)>\lz\})\le \dfrac
{C_p}{\lz^p}\dint_\rz|f(x)|^p\wz(x)dx.$$
\end{lem}

From Lemma 2.3, we know that $M_\vz$ may be not bounded on
$L^p(\wz)$ for all $\wz\in A_p(\vz)$ and $1<p<\fz$. We now need to
define a variant maximal operator $M_{\vz,\eta}$ for $0<\eta<\fz$
by
$$M_{\vz,\eta}f(x)=\dsup_{x\in Q}\dfrac 1{\vz(|Q|)^{\eta}|Q|}\dint_Q|f(y)|\,dy.$$
\begin{prop}\label{p2.1.}\hspace{-0.1cm}{\rm\bf 2.1.}\quad
Let $1<p<\fz$, $p'=p/(p-1)$ and suppose that $\wz\in A_p(\vz)$.
 There
exists a constant $C_p>0$ such that
$$\|M_{\vz,p'}f\|_{L^p(\wz)}\le C_p\|f\|_{L^p(\wz)}.$$
\end{prop}
The proof  can be found in \cite{tl} and  \cite {s}. We omit the
details here.

As a consequence of Proposition 2.1, we have
\begin{cor}\label{c2.1.}\hspace{-0.1cm}{\rm\bf 2.1.}\quad
Let $1\le p<\fz$ and $\wz\in A_p(\vz)$. Let $\Psi$ ba radial,
positive function with compact support and total integral 1. Set
$\Psi_t(x)=t^{-n}\Psi(x/t)$. Then
\begin{enumerate}
\item[(i)] $\dsup_{0<t<1}|f*\Psi_t(x)|\le C_\eta M_{\vz,\eta} f(x)$ for
$f\in L^p(\wz)$ and $0<\eta<\fz$;
\item[(ii)] $f*\Psi_t(x)\to f(x)$, as $t\to 0$, almost every for $f\in L^p(\wz)$;
\item[(iii)]$\|f*\Psi_t- f\|_{L^p(\wz)}\to 0$, as $t\to 0$, almost every for $f\in
L^p(\wz)$.
\end{enumerate}
\end{cor}

\begin{lem}\label{l2.4.}\hspace{-0.1cm}{\rm\bf 2.4.}\quad
If $\wz\in A_\fz(\vz)$,  then there exists positive constants
$C>0, \dz>0$ and $\dz_1>0$, such that for any
$Q=Q(x_0,r)\subset\rz$ with $r< 1$, we have
$$\l(\dfrac 1{|Q|}\dint_Q\wz^{1+\dz}dx\r)^{1/(1+\dz)}\le \dfrac
C{|Q|}\dint_Q\wz dx$$ and  for any measurable set $E\subset
Q=Q(x_0,r)$ with $r<1$, we have
$$\dfrac{\wz(E)}{\wz(Q)}\le C\l(\dfrac{|E|}{|Q|}\r)^{\dz_1}.$$
\end{lem}
Proof:\quad We first claim that for any $\az>0$, there exists a
positive constant  $c_0<1$ such that if $|A|/|Q|<\az$, then
$\frac{\wz(A)}{\wz(Q)}<c_0$ holds for any measurable set $A\subset
Q$.

Since $|A|<|Q|$, we have for any $x\in Q$,
$$M(\chi_{Q\setminus A})(x)>1-\az,$$
where $M$ denotes the standard Hardy-Littlewood maximal operator.

Note that if $r<1$, then there exists a constant $C>0$ such that
$M(\chi_{Q\setminus A})(x)\le CM_\vz(\chi_{Q\setminus A})(x)$ for
any $x\in Q$. From this and  by (iv) of Lemma 2.1, we have
$$M(\chi_{Q\setminus A})(x)\le C M_{\wz}(\chi^p_{Q\setminus
A})^{\frac 1p}(x).$$ By Lemma 2.2, we  get that
$$\begin{array}{cl}
\wz(Q)&\le\wz(\{x\in Q:\ M(\chi_{Q\setminus A})(x)>1-\az\})\\
&\le\wz(\{x\in Q:\ M_{\wz}(\chi_{Q\setminus A})(x)>C(1-\az)^p\})\\
&\le \dfrac C{(1-\az)^p}\dint_\rz\chi_{Q\setminus
A}(x)\wz(x)dx\\
&=\dfrac C{(1-\az)^p}(\wz(Q)-\wz(A)).
\end{array}$$
Obviously, $\frac C{(1-\az)^p}>1$. Hence, $\wz(A)\le c_0\wz(Q)$,
where $c_0<1$. Thus, the claim is proved. Using above claim and
adapting the standard proof in \cite{gr} and \cite{st}, we can
obtain the desired results. Thus, the proof of Lemma 2.4 is
finished.

Let $0<\eta<\fz$. We define the dyadic maximal operator
$M^\triangle_{\vz,\eta} f(x)$ by
$$M_{\vz,\eta}^\triangle f(x):=\dsup_{x\in   Q(dyadic\ cube)}\dfrac
1{\vz(|Q|)^\eta|Q|}\dint_Q|f(x)|\,dx.$$
\begin{lem}\label{l2.5.}\hspace{-0.1cm}{\rm\bf 2.5.}\quad
Let $f$ be a  locally integrable function on $\rz$, $\eta,\lz>0$,
and $\Omega_\lz=\{x\in \rz:\ M^\triangle_{\vz,\eta} f(x)>\lz\}$.
Then $\Omega_\lz$ may be written as a disjoint union of dyadic cubes
$\{Q_j\}$ with
\begin{enumerate}
\item[(i)] $\lz<(\vz(|Q_j|)^\eta|Q_j|)^{-1}\dint_{Q_j}|f(x)|\,dx,$
\item[(ii)] $(\vz(|Q_j|)^\eta|Q_j|)^{-1}\dint_{Q_j}|f(x)|\,dx\le
2^n\lz,$ for each cube $Q_j$. This has the immediate consequences:
\item[(iii)] $|f(x)|\le\lz \ {\rm for}\ a.e\ x\in
\rz\setminus\bigcup_jQ_j$
 \item[(iv)]$|\Omega_\lz|\le
\lz^{-1}\dint_\rz|f(x)|\,dx.$
\end{enumerate}
\end{lem}
The proof follows from the  argument of Lemma 1 in page 150 of
\cite{s}.

Let $0<\eta<\fz$. The dyadic sharp maximal operator
$M^\sharp_{\vz,\eta} f(x)$ is defined by
$$\begin{array}{cl}
M_{\vz,\eta}^{\sharp,\triangle} f(x)&:=\dsup_{x\in Q,r<1}\dfrac
1{|Q|}\dint_{Q(x_0,r)}|f(x)-f_Q|\,dx+ \dsup_{x\in Q,r\ge 1}\dfrac
1{\vz(|Q|)^\eta|Q|}\dint_{Q(x_0,r)}|f|\,dx\\
&\simeq\dsup_{x\in Q,r< 1}\dinf_C \frac
1{|Q|}\int_{Q(x_0,r)}|f(y)-C|\,dy+ \dsup_{x\in Q,r\ge1}\dfrac
1{\vz(|Q|)^\eta|Q|}\dint_{Q(x_0,r)}|f|\,dx
\end{array}$$
where $Q's$ denote dyadic cubes and $f_Q=\frac 1{|Q|}\int_Q
f(x)dx$. We define  sharp maximal operator $M^\sharp_{\vz,\eta}
f(x)$ as above if  dyadic cubes replaced by any cubes.

By Lemmas 2.4 and 2.5, we  establish the following ``good $\lz$"
inequality.
\begin{lem}\label{l2.6.}\hspace{-0.1cm}{\rm\bf 2.6.}\quad
Let $\wz\in A_\fz(\vz)$ and $0<\eta<\fz$. For a locally integrable
function $f$, and for $b$ and $\gz$ positive
$\gz<b<b_0=1/\vz(2^n)^\eta$, we have the following inequality
$$\wz(\{x\in \rz: M^\triangle_{\vz,\eta} f(x)>\lz, M^{\sharp,\triangle}_{\vz,\eta}f(x)\le\gz
\lz\})\le a^{\dz_1}\wz(\{x\in \rz:\  M^\triangle_{\vz,\eta}
f(x)>b\lz\})\eqno(2.3)$$ for all $\lz>0$, where $a=2^n\gz/(1-\frac
b{b_0})$ and $\dz_1>0$ depends only on $\wz$(see Lemma 2.4).
\end{lem}
Proof:\quad We may assume that the set $\{x:\
M^\triangle_{\vz,\eta} f(x)>b \lz\}$ has finite measure, otherwise
the inequality (2,3) is obvious. By Lemma 2.5, that this set is
the union of disjoint maximal cubes $\{Q_j\}$. We let
$Q=Q(x_0,r)=Q_j$ denote one of these cubes. Thus, we need only to
show that
$$\wz(\{x\in Q: M^\triangle_{\vz,\eta} f(x)>\lz, M^{\sharp,\triangle}_{\vz,\eta}f(x)\le\gz
\lz\})\le a^{\dz_1}\wz(Q).\eqno(2.4)$$ We consider two cases about
sidelength $r$, that is, $r< 1$ and  $r\ge 1$.

Case 1. When $r< 1$, let $\wt Q\supset Q$ be the parent of $Q$, by
the maximality of $Q$ we have $|f|_{\wt Q}\le b\lz\vz(|\wt Q|)\le
b\lz/b_0$. So far all $x\in Q$ for which $ M_{\vz,\eta}^\triangle
f(x)>\lz$, it follows that $M_{\vz,\eta}^\triangle (f\chi_
Q)](x)>\lz$, and also that $M_{\vz,\eta}^\triangle [(f-f_{\wt
Q})\chi_Q](x)>(1-b/b_0)\lz$. By the weak type (1,1) of
$M_{\vz,\eta}^\triangle$(see (iv) of Lemma 2.5), we have
$$\begin{array}{cl}
|\{x\in Q: M^\triangle_{\vz,\eta} f(x)>\lz,
M^{\sharp,\triangle}_{\vz,\eta}f(x)\le\gz \lz\}|&\le \dfrac
1{(1-b/b_0)\lz}\dint_Q|f-f_{\wt
Q}|dx\\
&\le \dfrac 1{(1-b/b_0)\lz}\dint_{\wt Q}|f-f_{\wt
Q}|dx\\
&\le \dfrac {|\wt Q|}{(1-b/b_0)\lz}\dinf_{x\in
Q}M_{\vz,\eta}^\sharp f(x)\\
&\le \dfrac {2^n\gz| Q|}{1-b/b_0},
\end{array}$$
if the set in question is not empty. Thus (2.4) is proved in the
case $r<1$ by Lemma 2.4.

 Case 2. When $r\ge 1$, note that
 $$b\lz<\dfrac 1{\vz(|Q|)^\eta|Q|}\dint_Q|f(y)|dy\le \dinf_{x\in Q}M_{\vz,\eta}^{\sharp,\triangle}
 f(x)\le\gz \lz,$$
but $\gz<b$, hence, the set in question is  empty. Thus (2.4) is
proved, and hence (2.3). the proof of Lemma 2.6 is complete.

 As a consequence of Lemma 2.6, we have the
following result.
\begin{cor}\label{c2.2.}\hspace{-0.1cm}{\rm\bf 2.2.}\quad
Let $1<p<\fz$, $\wz\in A_\fz(\vz)$ and $0<\eta<\fz$. Then there
exists a constant $C>0$ such that
$$\|M^\triangle_{\vz,\eta} f\|_{L^p(\wz)}\le
C\|M^{\sharp,\triangle}_{\vz,\eta} f\|_{L^p(\wz)}.$$
\end{cor}

Note that $|f(x)|\le M^\triangle_{\vz,\eta} f(x)\ a.e.\ x\in\rz$
and $M^{\sharp,\triangle}_{\vz,\eta} f(x)\le M^\sharp_{\vz,\eta}
f(x)$ for all $ x\in\rz$ for any $\eta>0$. By Corollary 2.2, we
have

\begin{prop}\label{p2.2.}\hspace{-0.1cm}{\rm\bf 2.2.}\quad
Let $1<p<\fz$, $\wz\in A_\fz(\vz)$, $0<\eta<\fz$ and $f\in
L^p(\wz)$, then
$$\| f\|_{L^p(\wz)}\le \|M^\triangle_{\vz,\eta} f\|_{L^p(\wz)}\le
C\|M^{\sharp}_{\vz,\eta} f\|_{L^p(\wz)}.$$ \end{prop} A variant of
dyadic maximal operator and dyadic sharp maximal operator
$$M^\triangle_{\dz,\vz,\eta}
f(x)=M^\triangle_{\vz,\eta}(|f|^\dz)^{1/\dz}(x)$$ and
$$M_{\dz,\vz,\eta}^{\sharp}
f(x)=M^{\sharp}_{\vz,\eta}(|f|^\dz)^{1/\dz}(x),$$ which will
become the main tool in our scheme. Proposition 2.2 and Lemma 2.6
imply immediately that
\begin{prop}\label{p2.3.}\hspace{-0.1cm}{\rm\bf 2.3.}\quad
Let $1<p<\fz$, $\wz\in A_\fz(\vz)$, $0<\eta<\fz$ and  $\dz>0$.
\begin{enumerate}
\item[(a)] Let $\vz: (0,\fz)\to(0,\fz)$ be doubling, that is,
$\vz(2a)\le C\vz(a)$ for $a>0$. Then, there exists a constant $C$
depending upon the $A_\fz$ condition of $\wz$ and doubling
condition of $\vz$ such that
$$\begin{array}{cl}
\dsup_{\lz>0}\vz(\lz)&\wz(\{y\in\rz:\ M^\triangle_{\dz,\vz,\eta}
f(y)>\lz\}) \\
&\le C\dsup_{\lz>0}\vz(\lz)\wz(\{y\in\rz:\
M^{\sharp}_{\dz,\vz,\eta}  f(y)>\lz\})\end{array}$$ for every
function such that the left hand side is finite. \item[(b)] If
$f\in L^p(\wz)$, then
$$\| f\|_{L^p(\wz)}\le\|M^\triangle_{\dz,\vz,\eta} f\|_{L^p(\wz)}\le
C\|M^{\sharp}_{\dz,\vz,\eta} f\|_{L^p(\wz)}.$$
\end{enumerate}\end{prop}

\begin{center} {\bf 3. The proof of Theorem 1.1}\end{center}
In this section, our proof follows from \cite{m} and \cite{p}. It
is worth  pointing out that  our proof here is rather complex,
which is of independent interest. As in \cite{m}, we first give a
result for any $ L^{-\fz}_{1,0}$ pseudo-differential operator.
\begin{lem}\label{l3.1}\hspace{-0.1cm}{\rm\bf3.1.}\quad
Suppose $A\in L^{-\fz}_{1,0}$ and $0<\eta<\fz$. Then exists a
constant $C>0$ such that for all $x_0\in\rz$ and all $f\in \cal
S$, $$M^{\sharp}_{\vz,\eta} (Af)(x_0)\le C M_{\vz,\eta} f(x_0).$$
\end{lem}
Proof:\quad If $a(x,\xi)$ is the symbol of $A$, then for any real
number $m$, and any multi-indices $\az$ and $\bz$,
$$|D_x^\az D_\xi^\bz a(x,\xi)|\le
C_{\az,\bz,m}(1+|\xi|)^{m-|\bz|}.$$ we can write the operator as
follows:
$$Af(x)=\dint_\rz\hat f(\xi)a(x,\xi)e^{2\pi
ix\cdot\xi}\,d\xi=\dint_\rz K(x,x-y)f(y)\,dy,$$ where
$$K(x,x-y)=\dint_\rz a(x,\xi)e^{2\pi
i(x-y)\cdot\xi}\,d\xi.$$ Obviously, for any $k>0$, there is a
constant $C_k>0$ such that
$$|K(x,x-y)|\le C_k(1+|x-y|)^{-k},\ {\rm for \ all}\ x\in\rz.$$
Then
$$\begin{array}{cl}
|Af(x)|&\le \dint_\rz|K(x,x-y)||f(y)|\,dy\\
&\le C_k\dint_\rz\dfrac {|f(y)|}{(1+|x-y|)^k}\,dy.
\end{array}$$
From this, we obtain that
$$\|Af\|_{L^1(\rz)}\le C\|f\|_{L^1(\rz)}.$$
For any cube $Q=Q(x_1,r)\ni x_0$, we write
$f:=f\chi_{2Q}+f\chi_{\rz\setminus 2Q}:=f_1+f_2$. Then we have
$$\begin{array}{cl}
\dfrac 1{\vz(|Q|)^\eta|Q|}\dint_Q|Af(x)|dx &\le \dfrac
1{\vz(|Q|)^\eta|Q|}\dint_Q|Af_1(x)|dx+\dfrac 1{\vz(|Q|)^\eta|Q|}\dint_Q|Af_2(x)|dx\\
&\le \dfrac
C{\vz(|Q|)^\eta|Q|}\dint_\rz|f_1(y)|\,dy\\
&\qquad\quad+\dfrac {C_k}{\vz(|Q|)^\eta|Q|}\dint_Q\dint_\rz\dfrac {|f_2(y)|}{(1+|x-y|)^k}\,dydx\\
&\le \dfrac
C{\vz(|Q|)^\eta|Q|}\dint_{2Q}|f(y)|\,dy\\
&\qquad\quad+\dfrac {C_k}{\vz(|Q|)^\eta|Q|}\dint_Q\dint_{|y-x_0|>r}\dfrac {|f(y)|}{(1+|x_0-y|)^k}\,dydx\\
&\le C M_{\vz,\eta} f(x_0)+C_k\dint_{r<|y-x_0|<1}\dfrac {|f(y)|}{(1+|x_0-y|)^k}\,dy\\
&\qquad\quad+C_k\dint_{1\le|y-x_0|}\dfrac
{|f(y)|}{(1+|x_0-y|)^k}\,dydx\\
&\le C M_{\vz,\eta} f(x_0).
\end{array}$$
Taking the supremum of the left side over all cubes $Q$ containing
$x_0$, we obtain the desired result.

Let us now turn to prove Theorem 1.1.

{\bf Proof of Theorem 1.1.}\quad To prove Theorem 1.1, from
propositions 2.1 and 2.3, we need only to show that for any
$p'<\eta<\fz$ and $0<\dz<1$ such that
$$M^\sharp_{\dz,\vz,\eta}(T f)(x_0)\le C_\eta M_{\vz,\eta}f(x_0)), \ {\rm a.e\ }\ x_0\in\rz,\eqno(3.1)$$
where $M^\sharp_{\dz,\vz,\eta}(T
f)(x)=M^\sharp_{\vz,\eta}(|Tf|^\dz)(x)^{1/\dz}.$

Fix $x_0\in\rz$ and let $x_0\in Q=Q(x_1,r)$. Decompose $f=f_1+f_2$,
where $f_1=f\chi_{\bar Q}$, where $\bar Q=Q(x_1,8r)$.

Case 1. When $r< 1$.
  Let $C_Q=|(T f_2)_Q|$.
Since $0<\dz<1$, then
$$\begin{array}{cl}
\l(\dfrac 1{|Q|}\dint_Q||T f(x)|^\dz-C_Q^\dz|\,dx\r)^{1/\dz}& \le
\l(\dfrac
1{|Q|}\dint_Q||T f(x)|-|(T f_2)_Q||^\dz\,dx\r)^{1/\dz}\\
&\le \l(\dfrac
1{|Q|}\dint_Q|T f(x)-(T f_2)_Q|^\dz\,dx\r)^{1/\dz}\\
&\le C\l(\dfrac
1{|Q|}\dint_Q|T f_1(x)|^\dz\,dx\r)^{1/\dz}\\
&\qquad+C\l(\dfrac
1{|Q|}\dint_Q|T f_2(x)-(T f_2)_Q|^\dz\,dx\r)^{1/\dz}\\
&\le C\l(\dfrac
1{|Q|}\dint_Q|T f_1(x)|^\dz\,dx\r)^{1/\dz}\\
&\qquad+C\dfrac
1{|Q|}\dint_Q|T f_2(x)-(T f_2)_Q|\,dx\\
&=I+II.
\end{array}$$
For I, we recall that $T$ is weak type $(1,1)$; see \cite{l}. Note
that
 $\vz(|\bar Q|)\sim
1$, by Kolmogorov's inequality(see\cite{p}), we then have
$$\begin{array}{cl}
I&\le \dfrac C{|Q|}\|T f_1\|_{L^{1,\fz}}\\
&\le \dfrac C{|\bar Q|}\dint_{\bar Q}|f(y)|\,dy\\
&\le CM_{\vz,\eta}f(x_0).
\end{array}$$
To deal with the second term, we shall also assume that $a(x,\xi)$,
the symbol of $T$, has compact $\xi$- support. The various constants
that occur in the following argument will not depend on the support
of $a$; at the end we show how to dispense with the assumption on
the support of $a$.

We decompose the operator $T$  into a sum of simpler operators. We
begin by fixing $\eta$, a nonnegative, radial, $C^\fz$ function of
compact support, defined in the $\xi-$space $\rz$, with the
properties that $\eta(\xi)=1$ for $|\xi|\le 1$ and $\eta(\xi)=0$ for
$|\xi|\ge 2$. Together with $\eta$, we define another function
$\phi$, by $\phi(\xi)=\eta(\xi)-\eta(2\xi)$. Then we have the
following  ``partitions of unity" of the $\xi-$space:
$$1=\eta(\xi)+\dsum_{j=1}^\fz\phi(2^{-j}\xi),\ {\rm all}\ \xi.$$
Now we can write
$$\begin{array}{cl}
Tf_2(x)&=\dint_\rz\hat f_2(\xi)a(x,\xi)e^{2\pi ix\cdot\xi}\,d\xi\\
&=\dint_\rz\hat f_2(\xi)a(x,\xi)\eta(\xi)e^{2\pi ix\cdot\xi}\,d\xi\\
&\qquad+\dsum_{j=1}^\fz\dint_\rz\hat f_2(\xi)a(x,\xi)\phi(2^{-j}\xi)e^{2\pi ix\cdot\xi}\,d\xi\\
&:=Af_2(x)+\dsum_{j=1}^\fz A_jf_2(x).
\end{array}$$
Obviously, $A\in L^{-\fz}_{1,0}$. Using  Lemma 3.1, we obtain that
$$\begin{array}{cl}
II&\le \dfrac C{|Q|}\dint_Q|A f_2(y)-(A f_2)_Q|\,dy+\dfrac
C{|Q|}\dint_Q|\dsum_{j=1}^\fz A_jf_2(x)-(\dsum_{j=1}^\fz A_jf_2)_Q|\,dx\\
&\le CM_{\vz,\eta}f(x_0)+C\dsum_{j=1}^\fz\dfrac 1{|Q|}\dint_Q|
A_jf_2(x)-( A_jf_2)_Q|\,dx.
\end{array}$$
It remains to examine the operators $A_j$,
$$A_jf_2(x)=\dint_\rz\dint_\rz a(x,\xi)\phi(2^{-j}\xi)e^{2\pi
i(x-y)\cdot\xi}\,d\xi\,dy.$$ The following lemma proved in \cite{m}
allows  to control the inner integral.
\begin{lem}\label{l3.2.}\hspace{-0.1cm}{\rm\bf 3.2.}\quad
Let $q(x,\xi)$ be a symbol of order $m$, and suppose $\phi\in
C_0^\fz(\rz)$ has support in $\{\xi:\ \frac12\le|\xi|\le 2\}$. If
$N\ge 0$, then there is a constant $C_N>0$ such that the
inequality
$$|y|^N\l|\dint_\rz q(x,\xi)\phi(2^{-j}\xi)e^{2\pi
ix\cdot\xi}\,d\xi\r|\le C_N 2^{j(n+m-N)}$$ holds for all $x$ and $y$
in $\rz$ and every integer $j\ge 1$.
\end{lem}

Let us continue to prove the theorem, noticing that
$$\begin{array}{cl}
\dfrac 1{|Q|}\dint_Q| A_jf_2(x)-( A_jf_2)_Q|\,dx&=\dfrac
1{|Q|}\dint_Q\l|\dfrac
1{|Q|}\dint_Q A_jf_2(x)- A_jf_2(z)dz\r|\,dx\\
&=\dfrac 1{|Q|}\dint_Q\l|\dfrac 1{|Q|}\dint_Q \dint_\rz
f_2(y)\dint_\rz\phi(2^{-j}\xi)\r.\\
&\quad\times\l.[a(x,\xi)e^{2\pi i(x-y)\cdot\xi}-a(z,\xi)e^{2\pi
i(z-y)\cdot\xi}\,d\xi\,dydz\r|dx.
\end{array}\eqno(3.2)$$
To estimate this last quantity, we consider two subcases:

Subcase 1. $2^jr\ge 1$. Taking $k_0$ such that $1/2<2^{k_0}r\le 1$.
Then (3.2) is dominated by
$$\begin{array}{cl}
2\dsum_{k=1}^\fz&\dfrac
1{|Q|}\dint_Q\dint_{2^kr\le|y-x_1|<2^{k+1}r}|f(y)|\l|\dint_\rz
\phi(2^{-j}\xi)a(x,\xi)e^{2\pi i(x-y)\cdot\xi}\,d\xi\r|\,dydx\\
&=2\dsum_{k=1}^{k_0}\dfrac
1{|Q|}\dint_Q\dint_{2^kr\le|y-x_1|<2^{k+1}r}|f(y)|\l|\dint_\rz
\phi(2^{-j}\xi)a(x,\xi)e^{2\pi i(x-y)\cdot\xi}\,d\xi\r|\,dydx\\
& +2\dsum_{k=k_0+1}^\fz\dfrac
1{|Q|}\dint_Q\dint_{2^kr\le|y-x_1|<2^{k+1}r}|f(y)|\l|\dint_\rz
\phi(2^{-j}\xi)a(x,\xi)e^{2\pi i(x-y)\cdot\xi}\,d\xi\r|\,dydx\\
&:=I_1+I_2.
\end{array}$$
For $I_1$, by Lemma 3.2 with $N=n+1$ and $m=0$, we obtain that
$$\begin{array}{cl}
I_1&\le C\dsum_{k=1}^{k_0}\dint_Q\dfrac
{2^{nk}}{|Q_k|}\dint_{2^kr\le|y-x_1|<2^{k+1}r}\dfrac{|f(y)|}{|x-y|^{n+1}}|x-y|^{n+1}\\
&\qquad\times\l|\dint_\rz
\phi(2^{-j}\xi)a(x,\xi)e^{2\pi i(x-y)\cdot\xi}\,d\xi\r|\,dydzdx\\
&\le C\dsum_{k=1}^{k_0} r^{-1}2^{-k}2^{-j}\dfrac
1{|Q_k|}\dint_{Q_k}|f(y)|dy\\
&\le Cr^{-1}2^{-j}M_{\vz, \eta}f(x_0).
\end{array}$$
For $I_2$, by Lemma 3.2 with $N=n+n\eta\az_0+1$ and $m=0$, we have
$$\begin{array}{cl}
I_2&\le C\dsum_{k=k_0+1}^\fz\dint_Q\dfrac
{2^{nk}}{|Q_k|}\dint_{2^kr\le|y-x_1|<2^{k+1}r}\dfrac{|f(y)|}{|x-y|^{n+n\eta\az_0+1}}|x-y|^{n+n\eta\az_0+1}\\
&\qquad\times\l|\dint_\rz
\phi(2^{-j}\xi)a(x,\xi)e^{2\pi i(x-y)\cdot\xi}\,d\xi\r|\,dydx\\
&\le C\dsum_{k=k_0+1}^\fz r^{-1}2^{-k}2^{-j-n\eta\az_0}\dfrac
1{\vz(|Q_k|)^\eta|Q_k|}\dint_{Q_k}|f(y)|dy\\
&\le Cr^{-1}2^{-j}M_{\vz, \eta}f(x_0).
\end{array}$$
Subcases 2. $2^jr<1$. We write
$$\begin{array}{cl}
a(x,\xi)&e^{2\pi i(x-y)\cdot\xi}-a(z,\xi)e^{2\pi i(z-y)\cdot\xi}\\
&=
\dsum_{l=1}^n(x_l-z_l)\dint_0^1\dfrac{\pz a}{\pz x_l}(x(t),\xi)
e^{2\pi i(x(t)-y)\cdot\xi}+2\pi i\xi_la(x(t),\xi)e^{2\pi
i(x(t)-y)\cdot\xi}dt,
\end{array}$$ where $x(t)=z+t(x-z)$.

Using this last expression and the facts:
\begin{enumerate}
\item[(a)] $\pz a/\pz x_l$ is a symbol of order $0$;
\item[(b)] $\xi_l a(x,\xi)$ is a symbol of order $1$;
\item[(c)]$|x_l-z_l|\le r$ since both $x$ and $z$ in $Q$; and
\item[(d)] if $2^kr\le|y-x_1|\le 2^{k+1}r$, then
$2^{k-1}r\le|x(t)-y|\le 2^{k+2}r$ since $x(t)\in Q$,
\end{enumerate}
we let $k_0$ as above, by Lemma 3.2 again,  then (3.2) is
dominated by
$$\begin{array}{cl}
&C\dsum_{k=1}^{k_0}\dfrac
1{|Q|}\dint_Q\dint_{2^kr\le|y-x_1|<2^{k+1}r}\dfrac{|f(y)|}{|x_1-y|^{n+1/2}}
\dsum_{l=1}^n|x_l-z_l|\dint_0^1|x(t)-y|^{n+1/2}\\
&\ \times\l|\dint_\rz \phi(2^{-j}\xi)\l[\dfrac{\pz a}{\pz
x_l}(x(t),\xi) e^{2\pi i(x(t)-y)\cdot\xi}+2\pi
i\xi_la(x(t),\xi)e^{2\pi i(x(t)-y)\cdot\xi}\r]d\xi\r|dtdydz\\
\end{array}$$
$$\begin{array}{cl}
&\ + \dsum_{k=k_0+1}^\fz\dfrac
C{|Q|}\dint_Q\dint_{2^kr\le|y-x_1|<2^{k+1}r}\dfrac{|f(y)|\dsum_{l=1}^n|x_l-z_l|}{|x_1-y|^{n+\frac12+n\eta\az_0}}
\dint_0^1|x(t)-y|^{n+\frac12+n\eta\az_0}\\
&\ \times\l|\dint_\rz \phi(2^{-j}\xi)\l[\dfrac{\pz a}{\pz
x_l}(x(t),\xi) e^{2\pi i(x(t)-y)\cdot\xi}+2\pi
i\xi_la(x(t),\xi)e^{2\pi i(x(t)-y)\cdot\xi}\r]d\xi\r|dtdydz\\
&\le C\dsum_{k=1}^{k_0}\dfrac
{2^{nk}}{|Q_k|}\dint_{Q_k}|f(y)|dyr^n(2^kr)^{-n-1/2}r(2^{-j/2}+2^{j/2})\\
&\ + C\dsum_{k=k_0+1}^\fz\dfrac
{2^{nk}}{\vz(|Q_k|)^\eta|Q_k|}\dint_{Q_k}|f(y)|dyr^n(2^kr)^{-n-1/2}r2^{j(-n\eta\az_0+1/2)}\\
&\le CM_{\vz,\eta}f(x_0)r^{1/2}2^{j/2}\dsum_{k=1}^\fz 2^{-k/2}\\
&\le Cr^{1/2}2^{j/2}M_{\vz,\eta}f(x_0).
\end{array}$$
Putting the two subcases together, we get that
$$\begin{array}{cl}
\dsum_{j=1}^\fz\dfrac 1{|Q|}\dint_Q| A_jf_2(x)-( A_jf_2)_Q|\,dx&\le
C\l(\dsum_{2^jr\ge1}r^{-1}2^{-j}+\dsum_{2^jr<1}r^{1/2}2^{j/2}\r)M_{\vz,\eta}f(x_0)\\
&\le CM_{\vz,\eta}f(x_0).
\end{array}$$
Case 2. When  $r\ge 1$, write $\rho_1:=\eta/\dz\ge\eta+1$, we have
 $$\begin{array}{cl}
 \l(\dfrac 1{\vz(|Q|)^{\eta}|Q|}\dint_Q|T
f(y)|^\dz\,dy\r)^{1/\dz}& \le \dfrac C{\vz(|Q|)^{\rho_1}}\l(\dfrac
1{|Q|}\dint_Q|T f_1(y)|^\dz\,dy\r)^{1/\dz}\\
&\qquad+\dfrac C{\vz(|Q|)^{\rho_1}}\l(\dfrac 1{|Q|}\dint_Q|T
f_2(y)|^\dz\,dy\r)^{1/\dz}\\
 &:=II_1+II_2.
\end{array}$$
For $II_1$, similar to $I_1$, we have
$$\begin{array}{cl}
I_1&\le \dfrac C{\vz(|Q|)^{\rho_1}}\dfrac 1{|Q|}\|T f_1\|_{L^{1,\fz}}\\
&\le \dfrac C{\vz(|Q|)^\eta| Q|}\dint_{\bar Q}|f(y)|\,d\mu(y)\\
&\le CM_{\vz,\eta}f(x).
\end{array}$$
Finally, for $II_2$, adapting the argument of  $I_2$ in subcase 1,
 note that $r\ge 1$, we obtain
$$\begin{array}{cl}
II_2 &\le C\dsum_{j=1}^\fz\dsum_{k=1}^\fz\dint_Q\dfrac
{2^{nk}}{|Q_k|}\dint_{2^kr\le|y-x_1|<2^{k+1}r}\dfrac{|f(y)|}{|x-y|^{n+\eta\az_0+1}}|x-y|^{n+n\eta\az_0+1}\\
&\qquad\qquad\times\l|\dint_\rz
\phi(2^{-j}\xi)a(x,\xi)e^{2\pi i(x-y)\cdot\xi}\,d\xi\r|\,dydx\\
\end{array}$$
$$\begin{array}{cl}
&\le C\dsum_{j=1}^\fz\dsum_{k=1}^\fz
r^{-1}2^{-k}2^{-j-n\eta\az_0}\dfrac
1{\vz(|Q_k|)^\eta|Q_k|}\dint_{Q_k}|f(y)|dy\\
&\le C\dsum_{j=1}^\fz r^{-1}2^{-j}M_{\vz, \eta}f(x_0)\\
&\le CM_{\vz, \eta}f(x_0).
\end{array}$$
Since $r\ge 1$.

 Finally, we show how to dispense with assumption on the support of $a$.
 Suppose now that the
assumption that $a(x,\xi)$, the symbol of $T$, has not compact
$\xi$-support. Let $b_j(x,\xi)$ be $a(x,\xi)$ multiplied by a
smooth cutoff function which is $1$ when $|\xi|\le 2^j$ and $0$
when $|\xi|\ge 2^{j+1}$. Let $B_j$ be the pseudo-differential
operator whose symbol is $b_j(x,\xi)$. Since $b_j(x,\xi)\to
a(x,\xi)$ as $j\to\fz$, the dominated convergence theorem implies
that $B_jf(x)\to Tf(x)$ for all $x$. Another application of the
dominated convergence theorem shows that for each cube $Q$ and
$0<\dz<1$,
$$\l(\dfrac 1{|Q|}\dint_Q||B_j f(x)|^\dz-|(B_jf_2)_Q|^\dz|\,dx\r)^{1/\dz}
\to \l(\dfrac 1{|Q|}\dint_Q||T
f(x)|^\dz-|(Tf_2)_Q|^\dz|\,dx\r)^{1/\dz}$$ and $$\l(\dfrac
1{\vz(|Q|)^{\eta}|Q|}\dint_Q|B_j f(y)|^\dz\,dy\r)^{1/\dz}\to
\l(\dfrac 1{\vz(|Q|)^{\eta}|Q|}\dint_Q|T f(y)|^\dz\,dy\r)^{1/\dz}.$$
Applying our previous result to the operators $B_j$, and taking the
limit as $j\to\fz$, we see that
$$\l(\dfrac 1{|Q|}\dint_Q||T f(x)|^\dz-|(Tf_2)_Q|^\dz|\,dx\r)^{1/\dz}\le CM_{\vz,
\eta}f(x_0),$$ if $|Q|<1$, and $$ \l(\dfrac
1{\vz(|Q|)^{\eta}|Q|}\dint_Q|T f(y)|^\dz\,dy\r)^{1/\dz}\le CM_{\vz,
\eta}f(x_0),$$ if $|Q|\ge 1$.

 Thus (3.1) holds. Theorem 1.1 is proved.

By Proposition 2.3 and (3.1), the weighted weak-type (1,1)
estimate for the operator $T$ is obtain as follows:
\begin{thm}\label{t3.1.}\hspace{-0.1cm}{\rm\bf 3.1.}\quad Let $T\in
L_{1,0}^0$ and  $\wz\in A_1(\vz)$. There exists a constant $C>0$
such that for any $\lz>0$,
$$\wz(\{x\in\rz: \ |Tf(x)|>\lz\})\le \dfrac
{C}\lz\dint_\rz |f(x)|\wz(x)d\mu(x).$$
\end{thm}

{\bf Proof of Corollary 1.1.}\quad  It is easy to see that
$\wz_1(x)=(1+|x|)^{\az_1}$ if $-n\az_0<\gz_1<n\az_0\in A_1(\vz)$
and $\wz_2(x)=|x|^{\az_2}\in A_1(\vz)$ for $-n< \az_2<n(p-1)$.
Then $\wz=w_2^\az w_1^{1-\az}\in A_p(\vz)$ for $1<p<\fz$ and
$0<\az<1$ by (iii) of Lemma 2.1, by Theorem 1.1, we have
$$\|Tf\|_{L^p(\wz)}\le C\|f\|_{L^p(\wz)}.$$
Form this and  by the arbitrary of $\az_0$, if
$\wz(x)=(1+|x|)^{\gz_1 }|x|^{-\gz_2}$ with $\gz_1\in\rr$ and $-n<
\gz_2<n(1-1/p)$, we then have
$$\|Tf\|_{L^p(\wz)}\le C\|f\|_{L^p(\wz)}.$$
Thus, Corollary 1.1 is proved.

\begin{center} {\bf 4. The proof of Theorem 1.2}\end{center}
In this section, our proof follows from \cite{st} and \cite{p},
which is  different from Section 3.

We first recall some basic definitions and facts about Orlicz
spaces, referring to \cite{r} for a complete account.

 A function $B(t): [0,\fz)\to [0,\fz)$ is called a Young function if
it is continuous, convex, increasing and satisfies $\Phi(0)=0$ and
$B\to \fz$ as $t\to\fz$. If $B$ is a Young function, we define the
$B$-average of a function $f$ over a cube $Q$ by means of the
following Luxemberg norm:
$$\|f\|_{B,Q}=\dinf\l\{\lz>0:\ \dfrac
1{|Q|}\dint_QB\l(\dfrac {|f(y)|}\lz\r)\,dy\le 1\r\}.$$ The
generalized H\"older's inequality
$$\dfrac
1{|Q|}\dint_B|fg|\,dy\le \|f\|_{B,Q}\|g\|_{\bar B,Q}$$ holds, where
$\bar B$ is the complementary Young function associated to $B$. And
we define the corresponding maximal function
$$M_B f(x)=\dsup_{Q:x\in Q}\|f\|_{B,Q}$$
and for $0<\eta<\fz$
$$M_{B,\vz,\eta} f(x)=\dsup_{Q:x\in Q}\vz(|Q|)^{-\eta}\|f\|_{B,Q}.$$

 The  example that we are going to use is
$B(t)=t(1+log^+t)$ with the maximal function denoted by $M_{LlogL}$.
The complementary Young function is given by $\bar B(t)\approx e^t$
with the corresponding maximal function denoted by $M_{exp L}$.

We need the following  several key lemmas.
\begin{lem}\label{l4.1.}\hspace{-0.1cm}{\rm\bf 4.1.}\quad
Let $ 0<\eta<\fz$ and $M_{\vz,\eta/2}f$ be locally integral. Then
there exists positive constants $C_1$ and $C_2$  independent of
$f$ and $x$ such that
$$ C_2M_{\vz,\eta} M_{\vz,\eta}f(x)\le M_{L\log L,\vz,\eta}f(x)\le
C_1M_{\vz,\eta/2} M_{\vz,\eta/2}f(x).$$
\end{lem}
Proof:\quad  It suffices to show that for any cube $Q\ni x$, there
is a constant $C>0$ such that
$$\vz(|Q|)^{-\eta}\|f\|_{L\log L,Q}\le C
\vz(|Q|)^{-\eta/2}|Q|^{-1}\dint_QM_{\vz,\eta/2}f(y)dy.$$ That is,
$$\|f\|_{L\log L,Q}\le C\vz(|Q|)^{\eta/2}{|Q|}^{-1}\dint_QM_{\vz,\eta/2}f(y)dy.\eqno(4.1)$$
In fact, by homogeneity we can take $f$ with $\|f\|_{L\log L,Q}=1$
which implies
$$\begin{array}{cl}
1&\le \dfrac C{|Q|}\dint_Q|f(y)|(1+\log^+(|f(y)|))dy\\
&\le \dfrac C{|Q|}\dint_Q|f(y)|\dint_1^{|f(y)|+1}\dfrac {dt}tdy\\
&\le \dfrac C{|Q|}\dint_1^\fz\dint_{\{x\in Q:
|f(y)|>t-1\}}|f\chi_Q(y)|dy\dfrac {dt}t\\
&\le \dfrac C{|Q|}\dint_0^\fz\dint_{\{x\in Q:
|f(y)|>t\}}|f\chi_Q(y)|dy\dfrac {dt}t\\
&\le \dfrac C{|Q|}\dint_0^\fz|\{x\in Q:\ M(f\chi_Q)(x)>t\}|dt\\
&=\dfrac C{|Q|}\dint_Q M(f\chi_Q)(x)dx\\
&\le C\dfrac {\vz(|Q|)^{\eta/2}}{|Q|}\dint_Q
M_{\vz,\eta/2}(f)(x)dx,
\end{array}$$
since $M(f\chi_Q)(x)\le C\vz(|Q|)^{\eta/2}M_{\vz,\eta/2}(f)(x)$
for all $x\in Q$. Thus, (4.1) is proved. Now let us turn to prove
$$M_{\vz,\eta}M_{\vz,\eta}f(x)\le CM_{L\log L,\vz,\eta}f(x).\eqno(4.2)$$  For any fixed $x\in\rz$ and any fixed
cube $Q\ni x$, write $f=f_1+f_2$, where $f_1=f\chi_{3Q}$. Thus,
$$\begin{array}{cl}
\vz(|Q|)^{-\eta}|Q|^{-1}\dint_Q|M_{\vz,\eta}f(y)|dy
&\le\vz(|Q|)^{-\eta}|Q|^{-1}\dint_Q|M_{\vz,\eta}f_1(y)|dy\\
&\quad+\vz(|Q|)^{-\eta}|Q|^{-1}\dint_Q|M_{\vz,\eta}f_2(y)|dy\\
&:=I+II.
\end{array}$$
For $I$, we know that for all $g$ with $\supp \ g\subset Q$ (see
\cite{p})
$$
\dfrac 1{|Q|}\dint_QMg(y)dy\le C\|g\|_{L\log L,Q}.$$ From this and
note that $M_{\vz,\eta}f(x)\le Mf(x)$, we get
$$\begin{array}{cl}
I&\le
C\vz(|Q|)^{-\eta}|3Q|^{-1}\dint_{3Q}|Mf_1(y)|dy\\
&\le C\vz(|Q|)^{-\eta}\|f\|_{L\log
L,3Q}\\
&\le CM_{L\log L,\vz,\eta}.\end{array}$$
 Next let us estimate $II$. It is easy to see that for all $y,\
 z\in Q$, we have
 $$M_{\vz,\eta}f_2(y)\le
 CM_{\vz,\eta}f(z).\eqno(4.3)$$
In fact,  for any cube $Q'\ni y$ and $Q'\bigcap(\rz\setminus
3Q)\not=\O$, noticing that $z\in Q\subset 3Q'$, we have
$$\begin{array}{cl}
\dfrac 1{\vz(|3Q'|)^{\eta}|Q'|}\dint_{Q'}|f_2(y)|dy &\le
C\dfrac 1{\vz(|3Q'|)^{\eta}|3Q'|}\dint_{3Q'}|f_2(y)|dy\\
&\le CM_{\vz,\eta}f(z).\end{array}$$ Hence, (4.3) holds.

Using  (4.3), we obtain
$$\begin{array}{cl}
II&\le C\vz(|Q|)^{-\eta}|Q|^{-1}\dint_Q|M_{\vz,\eta}f_2(y)|dy\\
&\le C|Q|^{-1}\dint_Q\dinf_{z\in Q}M_{\vz,\eta}f(z)dz\\
&\le CM_{\vz,\eta}f(x)\le CM_{L\log L,\vz,\eta}f(x).\end{array}$$
Thus, (4.2) is proved.

\begin{lem}\label{l4.2.}\hspace{-0.1cm}{\rm\bf 4.2.}\quad
Let $b\in BMO$ and $1\le\eta<\fz$. Let $0<\dz<\ez<1$, then
$$M^\sharp_{\dz,\vz,\eta}([b,T]f)(x_0)\le C \|b\|_{BMO}(M_{\ez,\vz, \eta}(T f)(x_0)
+M_{L\log L,\vz,\eta}(f)(x_0)), \ {\rm a.e\ }\
x_0\in\rz\eqno(4.4)$$ holds for all  $f\in C_0^\fz(\rz)$.
\end{lem}
Proof:\quad Observe that for any constant $\lz$
$$[b,T]f(x)=(b(x)-\lz)T f(x)-T((b-\lz)f)(x).$$
As above we fix $x\in\rz$ and let $x\in Q=Q(x_0,r)$. Decompose
$f=f_1+f_2$, where $f_1=f\chi_{\bar B}$, where $\bar Q=Q(x,8r)$. Let
$\lz$ be a constant and $C_Q$ a constant to be fixed along the
proof.

To prove (4.2), we consider two cases about $r$.

Case 1. When  $r< 1$. Since $0<\dz<1$,  we then have
$$\begin{array}{cl}
\l(\dfrac
1{|Q|}\dint_Q|\r.&\l.|[b,T]f(y)|^\dz-|C_Q|^\dz|\,dy\r)^{1/\dz}\\&
\le \l(\dfrac
1{|Q|}\dint_Q|[b,T]f(y)-C_Q|^\dz\,dy\r)^{1/\dz}\\
& \le \l(\dfrac
1{|Q|}\dint_Q|(b(y)-\lz)T f(y)-T((b-\lz)f)(y)-C_Q|^\dz\,dy\r)^{1/\dz}\\
&\le C\l(\dfrac
1{|Q|}\dint_Q|(b(y)-\lz)T f(y)|^\dz\,dy\r)^{1/\dz}\\
&\qquad+C\l(\dfrac
1{|Q|}\dint_Q|T((b-\lz)f_1)(y)|^\dz\,dy\r)^{1/\dz}\\
&\qquad+C\l(\dfrac
1{|Q|}\dint_Q|T((b-\lz)f_2)(y)-C_Q|^\dz\,dy\r)^{1/\dz}\\
&=I+II+III.
\end{array}$$
 To deal with $I$, we first fix
$\lz=b_{\bar Q},$ the average of $b$ on $\bar Q$. Then for any
$1<q<\ez/\dz$, by the John-Nirenberg inequality of $BMO$, we
obtain that
$$\begin{array}{cl}
I &\le C \l(\dfrac 1{|\bar Q|}\dint_{\bar Q}|b(y)-b_{\bar Q}|^{\dz
q'}\,dy\r)^{q'/\dz}\l(\dfrac
1{|Q|}\dint_Q|T f(y)|^{\dz q}\,dy\r)^{\dz q}\\
&\le C \|b\|_{BMO}M_{\ez,\vz,\eta}(Tf)(x_0),
\end{array}$$ where $1/q'+1/q=1$.

For II, we recall that $T$ is weak type $(1,1)$. By Kolmogorov's
inequality, we then have
$$\begin{array}{cl}
II&\le \dfrac C{|Q|}\|T f_1\|_{L^{1,\fz}}\le \dfrac C{|\bar Q|}\dint_{\bar Q}|f(y)|\,dy\\
&\le CM_{\vz,\eta}f(x_0).
\end{array}$$
Finally, to estimate $III$, from page 241-250 in \cite{st}, we can
express $T$ as
$$Tf(x)=\dint_\rz K(x,y)f(y)dy,$$
and the kernel $K(x,y)$ satisfies
$$|K(x,y)|\le C_N|x-y|^{-n-N}\ {\rm for}\ N\ge 0,\eqno(4.5)$$ and
$$|\pz_xK(x,y)|\le
C_N|x-y|^{-n-1-N}\ {\rm for}\ N\ge 0.\eqno(4.6)$$
Fixed the value of $C_Q$ by taking $C_Q=(T((b-b_{\bar Q})f_2))_Q$,
the average of $T((b-b_{\bar Q})f_2)$ on $Q$. Let
$b_{Q_k}=b_{Q(x_0,2^{k+1}r)}$. Taking $k_0$ such that
$1/2<2^{k_0}r\le 1$. Then,
$$\begin{array}{cl}
III&\le \dfrac C{|Q|}\dint_Q|T((b-b_{\bar Q})f_2)(y)-(T((b-b_{\bar Q})f_2))_Q|\,dy\\
&\le \dfrac C{|Q|^2}\dint_Q\dint_Q\dint_{\rz\setminus \bar
Q}|K(y,\wz)-K(z,\wz)||(b(\wz)-b_{\bar Q})f(\wz)|d\wz dzdy\\
&\le \dfrac C{|Q|^2}\dint_Q\dint_Q\dint_{|x_0-\wz|>2r
}|K(y,\wz)-K(z,\wz)||(b(\wz)-b_{\bar Q})f(\wz)|d\wz dzdy\\
&\le \dfrac
C{|Q|^2}\dint_Q\dint_Q\dsum_{k=2}^\fz\dint_{2^kr\le|x_0-\wz|<2^{k+1}r
}|K(y,\wz)-K(z,\wz)||(b(\wz)-b_{\bar Q})f(\wz)|d\wz dzdy\\
&\le \dfrac
C{|Q|^2}\dint_Q\dint_Q\dsum_{k=2}^{k_0}\dint_{2^kr\le|x_0-\wz|<2^{k+1}r
}|K(y,\wz)-K(z,\wz)||(b(\wz)-b_{\bar Q})f(\wz)|d\wz dzdy\\
& \quad + \dfrac
C{|Q|^2}\dint_Q\dint_Q\dsum_{k=k_0+1}^\fz\dint_{2^kr\le|x_0-\wz|<2^{k+1}r
}\\
&\qquad\qquad\times|K(y,\wz)-K(z,\wz)||(b(\wz)-b_{\bar Q})f(\wz)|d\wz dzdy\\
&:=III_1+III_2.
\end{array}$$
Taking $N=0$ in (4.6), we obtain that
$$\begin{array}{cl}
III_1&\le C\dsum_{k=1}^{k_0}\dfrac {2^{-
 k}}{|Q_k|}\dint_{Q_k}|b(\wz)-b_{\bar Q}|f(\wz)|d\wz\\
&\le C\dsum_{k=1}^{k_0}\dfrac {2^{-
 k}}{|Q_k|}\l[\dint_{Q_k}|b(\wz)-b_{Q_k}|f(\wz)|d\wz+|b_{Q_k}-b_{\bar Q}|\dint_{Q_k}||f(\wz)|d\wz\r]\\
&\le C\dsum_{k=1}^\fz  2^{- k}\|b\|_{BMO}M_{L\log
L,\vz,\eta}(f)(x_0)
+ C\|b\|_{BMO}M_{\vz,\eta}(f)(x_0)\dsum_{k=1}^\fz k2^{- k}\\
&\le C\|b\|_{BMO}M_{L\log L,\vz,\eta}(f)(x_0),
\end{array}$$
 in last inequality we have  used
that $|b_{\bar Q}-b_{Q_k}|\le Ck\|b\|_{BMO}$ and
$M_{\vz,\eta}(f)(x)\le M_{L\log L,\vz,\eta}(f)(x)$.

Taking $N=n\eta\az_0$ in (4.6), we get that
$$\begin{array}{cl}
III_2&\le C\dsum_{k=k_0+1}^\fz\dfrac {2^{-
 k}}{\vz(|Q_k|)^\eta|Q_k|}\dint_{Q_k}|b(\wz)-b_{\bar Q}|f(\wz)|d\wz\\
&\le C\dsum_{k=k_0+1}^\fz\dfrac {2^{-
 k}}{\vz(|Q_k|)^\eta|Q_k|}\l[\dint_{Q_k}|b(\wz)-b_{Q_k}|f(\wz)|d\wz+|b_{Q_k}-b_{\bar Q}|\dint_{Q_k}||f(\wz)|d\wz\r]\\
&\le C\dsum_{k=1}^\fz  2^{- k}\|b\|_{BMO}M_{L\log L,\vz,\eta}(f)(x)
+ C\|b\|_{BMO}M_{\vz,\eta}(f)(x)\dsum_{k=1}^\fz k2^{- k}\\
&\le C\|b\|_{BMO}M_{L\log L,\vz,\eta}(f)(x_0).
\end{array}$$
Case 2. When  $r\ge 1$.  Since $0<\dz< \ez<1$ and $\eta\ge 1$, let
 $\rho_2:=\eta/\dz\ge\eta$, we then have
$$\begin{array}{cl}
\l(\dfrac
1{\vz(|Q|)^\eta|Q|}\dint_Q\r.&\l.|[b,T]f(y)|^\dz\,dy\r)^{1/\dz}\\
 & =
\dfrac 1{\vz(|Q|)^{\rho_2}}\l(\dfrac
1{|Q|}\dint_Q|(b(y)-\lz)T f(y)-T((b-\lz)f)(y)|^\dz\,dy\r)^{1/\dz}\\
&\le C\dfrac 1{\vz(|Q|)^{\rho_2}}\l(\dfrac
1{|Q|}\dint_Q|(b(y)-\lz)T f(y)|^\dz\,dy\r)^{1/\dz}\\
&\qquad+C\dfrac 1{\vz(|Q|)^{\rho_2}}\l(\dfrac
1{|Q|}\dint_Q|T((b-\lz)f_1)(y)|^\dz\,dy\r)^{1/\dz}\\
&\qquad+C\dfrac 1{\vz(|Q|)^{\rho_2}}\l(\dfrac
1{|Q|}\dint_Q|T((b-\lz)f_2)(y)|^\dz\,dy\r)^{1/\dz}\\
&=I+II+III.
\end{array}$$
 To deal with $I$, we first fix
$\lz=b_{\bar Q},$ the average of $b$ on $\bar Q$. Then for any
$1<q<\ez/\dz$, by the John-Nirenberg inequality of $BMO$,  we then
have
$$\begin{array}{cl}
I &\le C \l(\dfrac 1{|\bar Q|}\dint_{\bar Q}|b(y)-b_{\bar Q}|^{\dz
r'}\,dy\r)^{1/(r'\dz)}\dfrac 1{\vz(|Q|)^{\rho_2}}\l(\dfrac
1{|Q|}\dint_Q|T f(y)|^{\dz q}\,dy\r)^{1/(\dz q)}\\
&\le C \|b\|_{BMO}M_{\vz,\eta}(T f)(x_0),
\end{array}$$ where $1/q'+1/q=1$.

For II, since $T$ is weak type $(1,1)$. By Kolmogorov's inequality,
we then have
$$\begin{array}{cl}
II&\le \dfrac C{\vz(|Q|)^{\rho_2}} \dfrac 1{|Q|}\|T f_1\|_{1,\fz}\\
&\le \dfrac C{\vz(|Q|)^{\rho_2}}\dfrac 1{|\bar Q|}\dint_{\bar
Q}|f(y)|\,dy\le CM_{\vz,\eta}f(x_0).
\end{array}$$
Finally, for III we first  fix the value of $C_Q$ by taking
$C_Q=(T((b-b_{\bar Q})f_2))_Q$, the average of $T((b-b_{\bar
Q})f_2)$ on $B$. Let $b_{Q_k}=b_{Q(x_0,2^{k+1}r)}$. Taking
$N=n\eta\az_0+1$ in (4.5), note that $r\ge 1$, we obtain that
$$\begin{array}{cl}
II&\le \dfrac C{|Q|}\dint_Q|T((b-b_{\bar Q})f_2)(y)|\,dy\\
&\le \dfrac C{|Q|^2}\dint_Q\dint_Q\dint_{\rz\setminus \bar
Q}|K(y,\wz)||(b(\wz)-b_{\bar Q})f(\wz)|d\wz dzdy\\
&\le C\dsum_{k=1}^\fz  2^{- k}r^{-1}\|b\|_{BMO}M_{L\log
L,\vz,\eta}(f)(x_0)
+ C\|b\|_{BMO}M_{\vz,\eta}(f)(x)\dsum_{k=1}^\fz k2^{- k}r^{-1}\\
&\le C\|b\|_{BMO}M_{L\log L,\vz,\eta}(f)(x_0).
\end{array}$$

From  these, we get (4.4). Hence the proof is
 finished.

\begin{lem}\label{l4.3.}\hspace{-0.1cm}{\rm\bf 4.3.}\quad
Let $\wz\in A_1(\vz)$ and $\eta\ge 2$. There exists a positive $C$
such that for any function $f$ and $\lz>0$
$$\wz(\{x\in\rz: M_{L\log L,\vz,\eta}(f)(x)>\lz\})\le
C\dint_\rz\Phi(|f(y)|/\lz)\wz(y)dy.\eqno(4.5)$$
\end{lem}
Proof:\quad Let $K$ be any compact subset in $\{x\in\rz: M_{L\log
L,\vz,\eta}(f)(x)>\lz\})$. For any $x\in K$, by a standard
covering lemma, it is possible to choose  cubes $Q_1,\cdots, Q_m$
with pairwise disjoint interiors such that $K\subset
\bigcup_{j=1}^m 3Q_j$ and with $\|f\|_{L\log L,\vz,Q_j}>\lz$,
$j=1,\cdots, m$. This implies
$$\vz(|Q_j|)^2|Q_j|\le \dint_{Q_j}\dfrac {|f(y)|}\lz\l(1+\log^+\l(\dfrac
{|f(y)|}\lz\r)\r)\,dy.$$ From this, by (vi) in  Lemma 2.1 with
$p=1$ and $E=Q$, we obtain that
$$\begin{array}{cl}
\wz(3Q_j)&\le C\vz(|Q_j|)\wz(Q_j)\\
&=C\vz(|Q_j|)^2|Q_j|
\dfrac{\wz(Q_j)}{\vz(|Q_j|)|Q_j|} \\
&\le C\dfrac{\wz(Q_j)}{\vz(|Q_j|)|Q_j|}\dint_{Q_j}\dfrac
{|f(y)|}\lz\l(1+\log^+\l(\dfrac {|f(y)|}\lz\r)\r)\,dy\\
&\le C\dinf_{Q_j}\wz(x) \dint_{Q_j}\dfrac
{|f(y)|}\lz\l(1+\log^+\l(\dfrac {|f(y)|}\lz\r)\r)\,dy\\
&\le C\dint_{Q_j}\dfrac {|f(y)|}\lz\l(1+\log^+\l(\dfrac
{|f(y)|}\lz\r)\r)\wz(y)\,dy.
\end{array}$$
Thus, (4.5) holds, hence, the proof is complete.

Using Proposition 2.3 and Lemmas 4.1, 4.2 and 4.3, adapting the
same argument of \cite{p}, we can prove the following result.
\begin{lem}\label{l4.4.}\hspace{-0.1cm}{\rm\bf 4.4.}\quad
Let $b\in BMO, \eta>2$ and $\wz\in A_1(\vz)$. Then there exists a
positive constant $C$ such that for any smooth function $f$ with
compact support
$$\begin{array}{cl}
\dsup_{t>0}&\dfrac 1{\Phi(1/t)}\wz(\{x\in\rz:\
|[T,b]f(x)|>t\})\\
&\le C\Phi(\|b\|_{BMO}) \dsup_{t>0}\dfrac
1{\Phi(1/t)}\wz(\{x\in\rz:\ M_{L\log
L,\vz,\eta}f(x)>t\}).\end{array}$$
\end{lem}

\bigskip

Let us turn to prove the main theorems in this section.

{\bf Proof of Theorem 1.2:}\quad By Lemmas 4.1, 4.2 and Theorem
 1.1, we have
$$\begin{array}{cl}
\|[T,b]\|_{L^p(\wz)}&\le C\|M_{\vz,\eta}(T
f)\|_{L^p(\wz)}+C\|M_{L\log L,\vz,\eta}(f)\|_{L^p(\wz)}\\
&\le C\|T
f\|_{L^p(\wz)}+C\|M_{\vz,\eta/2} M_{\vz,\eta/2}f\|_{L^p(\wz)}\\
&\le C\|f\|_{L^p(\wz)},
\end{array}$$
 taking $\eta=2p'$.

The weighted weak-type (1,1) estimate for the commutator is the
 following.
\begin{thm}\label{t4.1.}\hspace{-0.1cm}{\rm\bf 4.1.}\quad Let $b\in
BMO$ and $\wz\in A_1(\vz)$. There exists a constant $C>0$ such
that for any $\lz>0$,
$$\wz(\{x\in\rz: \ |[T,b]f(x)|>\lz\})\le C \dint_\rz\dfrac
{|f(x)|}\lz\l(1+\log^+\l(\frac {|f(x)|}\lz\r)\r)\wz(x)dx.$$
\end{thm}
{Proof:}\quad  Let $\Phi(t)=t\log (e+t)$. By homogeneity, we need
only to show that
$$\wz(\{x\in\rz: \ |[T,b]f(x)|>1\})\le C\l(\dint_\rz\Phi(
|f(x)|)\wz(x)dx\r).$$ Indeed, using Lemmas 4.3 and 4.4 and
adapting the same argument of \cite{p}, we can obtain
$$\begin{array}{cl}
\wz(\{x\in\rz: &\ |[T,b]f(x)|>1\})\\
&\le C \dsup_{t>0}\dfrac 1{\Phi(1/t)}\wz(\{x\in\rz:\
|[T,b]f(x)|>t\})\\
&\le C\Phi(\|b\|_{BMO}) \dsup_{t>0}\dfrac
1{\Phi(1/t)}\wz(\{x\in\rz:\
M_{L\log L,\vz,\eta}f(x)>t\})\\
&\le C\dint_\rz\Phi( |f(x)|)\wz(x)dx,\end{array}$$  taking
$\eta\ge 2$.

 Thus, the proof of
Theorem 4.1 is complete.

Similar to the proof of Corollary 1.1, we can obtain   Corollary
1.2. We omit the details here.

 LMAM, School of Mathematical  Science

 Peking University

 Beijing, 100871

 P. R. China

\bigskip

 E-mail address:  tanglin@math.pku.edu.cn


\begin{center} {\bf References}\end{center}
\begin{enumerate}
\vspace{-0.3cm}
\bibitem[1]{al}J. Alverez and J. Hounie,
Estimates for the kernel and continuity properties of
pseudo-differential operators, Ark Math. 28(1990), 1-22.
\vspace{-0.3cm}
\bibitem[2]{at}
P. Auscher and M. Tayler, Paradifferential operators and
commutators estimates, Comm PDE. 20(1995), 1743-1775.
\vspace{-0.3cm}
\bibitem[3]{ct} S. Chanillo and A. Torchinsky, Sharp function and
weighted $L^p$ estimates for a class of pseudo-differential
operators, Ark Mat. 28(1986), 1-25.
\vspace{-0.3cm}
\bibitem[4]{cf} R. Coifman and C. Fefferman,
Weighted norm inequalities for maximal functions and singular
integrals, Studia Math. 51(1974), 241-250.
 \vspace{-0.3cm}
\bibitem[5]{gr} J. Garc\'ia-Cuerva and J. Rubio de Francia,
Weighted norm inequalities and related topics, Amsterdam- New York,
North-Holland, 1985. \vspace{-0.3cm}
\bibitem[6]{jn} F. John and L. Nirenberg,
On functions of bounded mean oscillation, Comm. Pure Appl. Math.
4(1961), 415-426. \vspace{-0.3cm}
\bibitem[7]{l} A. Laptev,
Spectral asympotics of a class of Fourier integral operators
(Russian), Trudy Moskov. Mat. Obshch. 43(1981), 92-115.
\vspace{-0.3cm}
\bibitem[8]{m}  N. Miller,
Weighted sobolev spaces and pseudodifferential operators with smoo
th symbols, Trans. Amer. Math. Soc. 269(1982), 91-109.
 \vspace{-0.3cm}
\bibitem[9]{p}C. P\'erez,
Endpoint estimates for commutators of singular integral operators,
J. Funct. Anal. 128(1995), 163-185. \vspace{-0.3cm}
\bibitem[10]{r} M.  Rao and Z.  Ren,
Theory of Orlicz spaces, Monogr. Textbooks Pure Appl. Math. 146,
Marcel Dekker, Inc., New York, 1991.
 \vspace{-0.3cm}
\bibitem[11]{s} E. Sawyer,
A characterization of a two-weight inequality for maximal operators,
Studia Math. 75(1982), 1-11.
\vspace{-0.3cm}
\bibitem[12]{st}  E.  Stein,
Harmonic Analysis: Real-variable Methods, Orthogonality, and
Oscillatory integrals. Princeton Univ Press. Princeton, N. J. 1993.
\vspace{-0.3cm}
\bibitem[13]{t} M. Taylor,
Pseudodifferential operators and nonlinear PDE. Boston:
Birkhau-ser, 1991. \vspace{-0.3cm}
\bibitem[14]{tl}
L. Tang,
 New  $BMO$  and weight functions associated with Schr\"odinger
 operators, preprint.
\end{enumerate}
\end{document}